\documentclass[11pt]{article}
\usepackage{latexsym}
\usepackage{amsmath}
\usepackage{amssymb}
\newtheorem{theorem}{Theorem}[section]

\newtheorem{proposition}[theorem]{Proposition}
\newtheorem{definition}[theorem]{Definition}
\newtheorem{lemma}[theorem]{Lemma}

\def\p{{\bf Proof.} \quad}
\def\q{\hfill\rule{1ex}{1ex}}

\begin{document}
\title{\bf Two coloring problems on matrix graphs}
\author{ { {\small\bf Zhe Han}
\thanks{email: hanzhe101@gmail.com.}
\quad {\small\bf Mei
Lu} \thanks{Correspondent author, email: mlu@math.tsinghua.edu.cn.}}\\
{\small Department of Mathematical Sciences, Tsinghua
University, Beijing 100084, China} }

\date{}

\maketitle

\baselineskip 16.3pt

\begin{abstract} In this paper,  we propose a new family of graphs,  matrix graphs,  whose vertex set $\mathbb{F}^{N\times n}_q$ is the set of all $N\times n$ matrices over a finite field $\mathbb{F}_q$ for any positive integers $N$ and $n$. And any two matrices share an edge if the rank of their difference is $1$.    Next,  we give some basic properties of such graphs and also consider two coloring problems on them.  Let $\chi'_d(N\times n, q)$ (resp. $\chi_d(N\times n, q)$) denote the minimum number of colors necessary to color the above matrix graph  so that no two vertices that are at a distance at most $d$ (resp. exactly $d$) get the same color. These two problems were proposed in the study of scalability of optical networks. In this paper, we  determine the exact value of $\chi'_d(N\times n,q)$ and give some upper and lower bounds on   $\chi_d(N\times n,q)$.

\end{abstract}

\vskip.3cm


\section{Introduction}

Let $q=p^m$, where $p$ is a prime. For any two positive integers $N$ and $n$,  
let $\mathbb{F}^{N\times n}_q$ be the set of all $N\times n$  matrices over  a finite field $\mathbb{F}_q$. In this paper we assume $n\leq N$ since the transpose of an $N\times n$ matrix becomes an $n\times N$ matrix.  The rank, or rank weight,  of a matrix $M\in \mathbb{F}^{N\times n}_q$, denoted by Rk$(M)$,  is defined as the algebraic rank  of this matrix over $\mathbb{F}_q$.
The rank  distance $d_{R}(M_1,M_2)$ between
$M_1$ and $M_2$ is the rank of their difference $M_1-M_2$.

Let $G=(V ,E)$ be a simple undirected graph with the vertex set $V$ and the edge set $E$. For two vertices $v_i$ and $v_j$ ($i\not = j$) in $V$, the distance between $v_i$ and $v_j$, denoted by $d(v_i,v_j)$, is the number of edges in a shortest path joining $v_i$ and $v_j$ and the diameter of $G$ is
the maximum distance between any two vertices of $G$.

 \begin{definition}
 A  matrix graph, denoted by $\mathcal{M}_{N\times n}(q)$, is defined to be an undirected graph with the vertex set $\mathbb{F}^{N\times n}_q$ and the edge set
 $$E_{N\times n}(q)=\{(M_1,M_2)|M_1, M_2\in\mathbb{F}^{N\times n}_q,  d_{R}(M_1,M_2)=1\}.$$
 \end{definition}

 Recently, hypercubes were extensively studied due to their versatile and efficient topological structures of interconnection networks in  \cite{Bhuyan,Fu,hanlu, Jamison, Kim, Klotz, Ngo,NDG,Ostergard,Skupie,Sullivan,Wan}.
 The so-called $q$-ary $n$-cube is an undirected graph with the vertex set, $n$-dimensional vector space over $\mathbb{F}_q$,
and the edge set, which contains all pairs such that their Hamming distance is $1$.  The Hamming distance is  the number of distinct coordinates between the pair.

Matrix graphs may be  better  than hypercubes and so it is worth studying their properties.
We first give a basic lemma which will be used to prove the following Proposition 1.3. It is an old result with
the first derivation of the formula in [14, p.455]; see also Lemma 2.1 in \cite{ling}.

\begin{lemma} \label{lem-3.1}
The number of $N \times n$ matrices  over $\mathbb{F}_q$ with rank $k$ is
 $ \frac{\prod_{i=0}^{k-1}(q^N-q^i)(q^n-q^i)}{\prod_{i=0}^{k-1}(q^k-q^i)}$.
\end{lemma}

\begin{proposition}
The matrix graph $\mathcal{M}_{N\times n}(q)$ is $\frac{(q^N-1)(q^n-1)}{q-1}$-regular connected graph of order $q^{Nn}$.
\end{proposition}

\begin{proposition}
The distance of any two distinct matrices $M_1$ and $M_2$ in $\mathbb{F}^{N\times n}_q$ has $d(M_1,M_2)=d_R(M_1,M_2)=\mbox{Rk}(M_1-M_2)$.
\end{proposition}

\begin{proposition}
The  diameter of $\mathcal{M}_{N\times n}(q)$ is $n$.
\end{proposition}

Note that $\mathcal{M}_{N\times n}(q)$ is not a bipartite graph and $\mathcal{M}_{N\times 1}(q)$ is a complete graph.

\begin{proposition}
The $\mathcal{M}_{N\times n}(q)$ is vertex transitive and then the edge connectivity is $\frac{(q^N-1)(q^n-1)}{q-1}$.
\end{proposition}

\noindent \p If $M$ is a fixed matrix in $\mathbb{F}^{N\times n}_q$, then the mapping
$$\rho_M~:~A\mapsto A+M$$
 is a permutation of the vertices of $\mathcal{M}_{N\times n}(q)$, where $+$ is the matrix addition. This mapping is an automorphism because for any two distinct matrices $M_1$ and $M_2$, $d_R( M_1,M_2)=1$ if and only if $d_R(M_1+M,M_2+M)=1$. There are $q^{Nn}$ such permutations and they form a subgroup of the automorphism group of $\mathcal{M}_{N\times n}(q)$. This subgroup acts transitively on $\mathbb{F}^{N\times n}_{q}$ because for any two vertices $M_1$ and $M_2$, the automorphism $\rho_{M_2-M_1}$ maps $M_1$ to $M_2$.\q

A coloring of $\mathcal{M}_{N\times n}(q)$ with $L$ colors is a mapping $\Gamma$ from the vertex set $\mathbb{F}^{N\times n}_q$ to ${\cal L}=\{1,2,\ldots,L\}$.
A $d$-distance (resp. exactly $d$-distance) coloring of $\mathcal{M}_{N\times n}(q)$ is to color the vertices of $\mathcal{M}_{N\times n}(q)$ such that
any two vertices with rank distance at most $d$ (resp. exactly $d$) have different colors. Note that
for a coloring of $\mathcal{M}_{N\times n}(q)$ with $L$ colors
\[\Gamma: \mathbb{F}^{N\times n}_q \longrightarrow \mathcal{L},\]
 it is a $d$-distance coloring of $\mathcal{M}_{N\times n}(q)$ if and only if for any two distinct vertices $M_1,M_2\in \mathbb{F}^{N\times n}_q$,
\[\Gamma(M_1)\neq \Gamma(M_2) \mbox{ if } d_R(M_1,M_2)\leq d\]
\noindent and it is an exactly $d$-distance coloring of $\mathbb{F}^{N\times n}_q$ if and only if for any two distinct vertices $M_1,M_2\in \mathbb{F}^{N\times n}_q$,
\[\Gamma(M_1)\neq \Gamma(M_2) \mbox{ if } d_R(M_1,M_2)= d.\]
Denote $\chi'_d(N\times n,q)$ (resp. $\chi_d(N\times n,q)$) as the minimum number of colors needed for a $d$-distance (resp. an exactly
$d$-distance) coloring of $\mathcal{M}_{N\times n}(q)$. Clearly, $\chi_{d}(N\times n,q)\leq \chi'_{d}(N\times n,q)$ for any $n$ and $N$.


These two coloring problems originally arose in the study of the scalability of optical networks \cite{Pavan}.  In the rest of our paper, we determine exact value of  $\chi'_{d}(N\times n,q)$ and give some bounds for $\chi_{d}(N\times n,q)$ for any $n$, $N$ and $d$.  

The $d$-distance coloring and exactly $d$-distance coloring of $\mathcal{M}_{N\times n}(q)$ are equivalent to certain partitions of $\mathbb{F}^{N\times n}_q$, which are related
to rank codes in coding theory.  Therefore,  we introduce rank codes  in the next section.

\section{Rank codes }

Algebraic coding theory can be considered as the theory of subsets of a  finite-dimensional space over a finite field equipped with a norm function.  The most known norm in coding theory is the Hamming weight of a vector. It turns out that the rank function of a matrix over a finite field can also be considered as a norm function. The interest in these codes is a consequence of their application in random network coding \cite{Koetter}. Explicitly, the concept of the rank metric was introduced by Loo-Keng Hua \cite{hua} as "Arithmetic distance".  Delsarte \cite{del} defined the rank distance on the set of bilinear forms (equivalently, on the set of rectangular matrices) and proposed the construction of optimal codes in bilinear form representation.  Gabidulin
\cite{gab} introduced the rank distance for the vector representation over extension fields and found connections between
rank codes in the vector representation and in the matrix representation.

In the matrix representation, rank codes are defined as subsets of a normed space $\{\mathbb{F}^{N\times n}_q$, $Rk\}$ of $N\times n$ matrices over
 $\mathbb{F}_q$.  The definitions of  norm and distance of any two matrices  are defined as above.
The rank distance of a rank code $\mathcal{M}\subset \mathbb{F}^{N\times n}_q$ is defined as the minimal pairwise distance: $d(\mathcal{M})=d=\min\{Rk(M_i-M_j): M_i,M_j\in \mathcal{M}, i\neq j\}$.

 The size $|\mathcal{M}|$ of related code with code distance $d$ satisfies  the Singleton bound
$|\mathcal{M}|\leq q^{N(n-d+1)}$.  Codes reaching this bound are called maximum rank codes, or, MRD codes.
A rank code $\mathcal{M}$ is called $\mathbb{F}_q$-linear
if $\mathcal{M}$ is a subspace of $\mathbb{F}^{N\times n}_q$.

In \cite{del} the construction of optimal codes is proposed.  Therefore for any $n\leq N$,  a linear rank code $\mathcal{M}$ with code distance $d$ reaches the singleton bound $|\mathcal{M}|=q^{N(n-d+1)}$. 

In the vector representation, rank codes are defined as subsets of a normed $n$-dimensional space
$\{\mathbb{F}^{n}_{q^N}, Rk\}$ of $n$-vectors over an extension field $\mathbb{F}_{q^N}$, where
the norm of a vector $\mathbf{v}\in \mathbb{F}^{n}_{q^N} $  is defined to be the column rank
$Rk(\mathbf{v}|\mathbb{F}_q)$ of this vector over $\mathbb{F}_q$, that is, the maximal number of coordinates of $\mathbf{v}$ which are linearly independent over the base field $\mathbb{F}_q$.  The  rank distance between two vectors $\mathbf{v}_1, \mathbf{v}_2$ is the column rank of their difference $Rk(\mathbf{v}_1-\mathbf{v}_2|\mathbb{F}_q)$.  The rank distance of a vector rank code $\mathcal{V}\subset \mathbb{F}^{n}_{q^N}$ is defined as the minimal pairwise distance:
$d(\mathcal{V})=d=\min\{Rk(\mathbf{v}_i-\mathbf{v}_j|\mathbb{F}_q): \mathbf{v}_i,\mathbf{v}_j\in \mathcal{V},i\neq j\}$.

A rank code $\mathcal{V}$ in the vector representation is called $\mathbb{F}_{q^N}$-linear if $\mathcal{V}$ is a subspace of $\mathbb{F}_{q^N}^{n}$.  Denote by $(N\times n,k,d)$ a code $\mathcal{V}$ over $\mathbb{F}_{q^N}$ of dimension $k\leq n$ and rank distance $d$.  Such a code can be described in terms of a full rank generator matrix $G_k$ over the extension field $\mathbb{F}_{q^N}$ of size $k\times n$. Code vectors $\{\mathbf{v}\}$ are all linear combinations of rows of this matrix.  Thus the size of a code is equal to $|\mathcal{V}|=q^{Nk}$.

Equivalently, a rank code $\mathcal{V}$ can be described in terms of a full rank parity-check matrix $H_{n-k}$ over $\mathbb{F}_{q^{N}}$ of size $(n-k)\times n$.   It satisfies the condition  $G_kH^{T}_{n-k}=O$, where $O$ is the all zero $k\times (n-k)$ matrix. Code vectors $(\mathbf{v})$ are all solutions of the linear system of equation $\mathbf{v}H^{T}_{n-k}=\mathbf{0}$.  We give the similar property of linear rank codes to the linear hamming codes as follows.

\begin{lemma}\label{lem:3.3}
Let $\mathcal{V}$ be a $(N\times n,k,d)$ linear rank code over $\mathbb{F}_{q^N}$ and  $H$ be a parity-check matrix for $\mathcal{V}$. Then the following statements are equivalent:

(i) $\mathcal{V}$ has rank distance $d$.

(ii) any $d-1$ columns of $H$ are linearly independent over $\mathbb{F}_{q^N}$ and $H$ has a $d$ columns
 that are linearly dependent over $\mathbb{F}_{q^N}$.
\end{lemma}

\noindent \p Let  $\mathbf{v}=(v_1,\ldots,v_n)\in \mathcal{V}$ be a codeword of rank weight $e$.  Thus there
are $n\times n$ invertible matrix $P$ such that $\mathbf{v}P=\mathbf{u}$, where $\mathbf{u}=(u_1,\ldots,u_e,0,\ldots,0)$ and where $u_i\neq 0$. So
$\mathbf{v}=\mathbf{u}P^{-1}$.     Then $\mathcal{V}$ contains a nonzero codeword
$\mathbf{v}$ of rank weight $e$ if and only if
\[\mathbf{0}=\mathbf{v}H^T=\mathbf{u}(H(P^{-1})^T)^T=u_1\mathbf{c}_1+\ldots+u_e\mathbf{c}_e,\]
which is true if and only if there are $e$ columns of $ H(P^{-1})^T$ (namely, $\mathbf{c}_{1}$,$\ldots$, $\mathbf{c}_{e}$)
that are linear dependent over $\mathbb{F}_{q^N}$ if and only if there are $e$ columns of $H$ that are linear dependent over $\mathbb{F}_{q^N}$.

To say that the rank distance of $\mathcal{V}$ is $\geq d$ is equivalent to saying that $\mathcal{V}$ does not contain any nonzero word of $\leq d-1$, which is in turn equivalent to saying that any $d-1$ columns of $H$ are linearly independent over $\mathbb{F}_{q^N}$.

Similarly,  to say that the rank distance of $\mathcal{V}$ is $\leq d$ is equivalent to saying that $\mathcal{V}$ contains a nonzero word of weight $\leq d$, which is in turn equivalent to saying that $H$ has $\leq d$ columns that are linear dependent over $\mathbb{F}_{q^N}$.

 The above discussion easily leads to the desired results.
\q

Likewise, general constructions of MRD codes in terms of parity-check matrices can be obtained \cite{gab}.

\begin{proposition}\label{prop:4}
Let $h_1,h_2,\ldots,h_n$ be a set of elements from the extension field $\mathbb{F}_{q^N}$
linearly independent over the base field $\mathbb{F}_q$. Let $s$ be a positive integer such that $\gcd(s,N)=1$. Then a parity matrix of the form
\[ H_{d-1}=\left(\begin{array}{cccc}
                        h_1& h_2& \ldots &h_n\\
                        h^{q^s}_{1} & h^{q^s}_{2} & \ldots & h^{q^s}_{n} \\
                        \ldots & \ldots & \ldots & \ldots \\
                        h^{q^{s(d-2)}}_{1} & h^{q^{s(d-2)}}_{2} & \ldots & h^{q^{s(d-2)}}_{n}
                        \end{array}\right )\]
defines an MRD $(N\times n,k,d)$ code with code length $n\leq N$, dimension $k=n-d+1$, and
rank distance $d=n-k+1$.
\end{proposition}

Equidistant codes in rank norm are related to the exactly $d$-distance coloring problem. A code is said to be equidistant if the distance between any distinct codewords is the same (say $d$). A code is called a constant-weight code if each non-zero codeword is of the same weight.  In \cite{Tuvi} the authors constructed an equidistant constant rank code as follows.

\begin{proposition} \label{prop:5}
There exists an equidistant constant rank code over $\mathbb{F}_q$ with
matrices of size ${n \choose 2}\times n$, rank $n-1$, rank distance $n-1$, and size $q^n-1$.
\end{proposition}

\section{Exact value of  $\chi'_{d}(N\times n,q)$ }

 A code $C$ over $\mathbb{F}_{q^N}$ of length $n$ and minimum rank distance at least $d$ is called an $(N\times n,\ge d)_q$ code. Let $A_q(N\times n,d)$ (resp. $A_q(N\times n,\ge d)$) denote the maximum size of an $(N\times n, d)_q$ (resp. $(N\times n,\ge d)_q$) code. A code $C$ is called an $[N\times n, k]_q$ linear code if $C$ is a $k$-dimensional subspace of $\mathbb{F}^{n}_{q^N}$.
 An $[N\times n,k]_q$ linear code with minimum distance $d$ is denoted by an  $[N\times n,k,d]_q$. For the fixed $N$, $n$, and minimum distance $d$, let $k(N\times n,d)_q$ denote the maximum dimension of an $[N\times n,k,d]_q$ code. A  code $C$ of length $n$ is called an $(N\times n, \overline{\{d\}})_q$ forbidden distance code if $d_R(\mathbf{u},\mathbf{v})\not= d$ for any two distinct codewords $\mathbf{u},\mathbf{v}\in C$. Given $N, n$, and $d$, let $Q(N\times n, d)_q$ denote the maximum size of an $(N\times n, \overline{\{d\}})_q$ forbidden rank distance code.

An $(N\times n,L,d)_q$-partition (resp. $(N\times n, L,\overline{\{d\}})_q$-partition) of $\mathbb{F}^{N\times n}_q$ is a set of subsets $\{\mathcal{M}_i\}_{i=1}^{L}$ of $\mathbb{F}^{N\times n}_q$ satisfying (i) $\mathcal{M}_i\bigcap \mathcal{M}_j=\emptyset$ for $i\neq j$ and $\bigcup_{i=1}^{L}\mathcal{M}_i=\mathbb{F}^{N\times n}_q$, (ii) each $\mathcal{M}_i$ is an $(N\times  n,\geq d)_q$ rank code (resp. $(N\times  n,\overline{\{d\}})_q$ forbidden distance rank code). It is well known that a $d$-distance coloring (resp. an exactly $d$-distance) of $\mathbb{F}^{N\times n}_q$ with $L$ colors is equivalent to an $(N\times n,L,d+1)_q$-partition (resp. $(N\times n,L,\overline{\{d\}})_q$-partition) of $\mathbb{F}^{N\times n}_q$. Hence $\chi'_{d}(N\times n,q)$ (resp. $\chi_{d}(N\times n,q)$) is the minimum number $L$ of subsets in any $(N\times n, L, d+1)_q$ (resp. $(N\times n, L,\overline{\{d\}})_q$-partition) of $\mathbb{F}^{N\times n}_q$.

Since $q^{Nn}=|\mathbb{F}^{N\times n}_q|=\sum_{i=1}^{L}|\mathcal{M}_i|\le L A_{q}(N\times n,\geq d+1)\leq L A_{q}(N\times n,d+1) $ and $A_q(N\times  n,d)$ is decreasing in $d$, we have
\begin{equation}\label{eq:basic}
 \chi'_{d}(N\times n,q)\geq \frac{q^{Nn}}{A_{q}(N\times n,d+1)}.
 \end{equation}

Note that for an $[N\times n,k,d+1]_q$ linear rank code $\mathcal{M}$, the cosets of $\mathcal{M}$ form an $(N\times n,q^{n-k},d+1)_q$-partition of $\mathbb{F}^{N\times n}_q$. Hence,
if there exists an $[N\times n,k,d+1]_q$ linear rank code, then
\begin{equation} \label{eq:code}
\chi'_d(N\times n,q)\leq q^{N(n-k)}.
 \end{equation}
In particular,
\begin{equation} \label{eq:max}
\chi'_d(N,n,q)\leq q^{N(n-k(n,d+1))}.
 \end{equation}
Furthermore, if $A_q(N\times n,d+1)=q^{N\times k(n,d+1)}$, i.e., $A_q(N\times n,d+1)$ is attained by an  $[N\times n,k,d+1]_q$ linear rank code, then
\begin{equation}\label{eq:final}
\chi'_{d}(N\times n,q)=q^{N(n-k(n,d+1))}.
\end{equation}

Summing the above discussion, we get our result on $\chi'_{d}(N\times n,q)$.
It is clear that $\chi'_{d}(N\times n,q)=1$ if $d\geq n$. 

\begin{theorem} For any positive integers $d\leq n$,  we have
$\chi'_{d}(N\times n,q)=q^{Nd}$.
\end{theorem}
\p
By Proposition \ref{prop:4}, we have $k(n,d+1)=n-d$ and so the required result is obtained by  Eq(\ref{eq:final}). 
\q

\section{Bounds of $\chi_{d}(N\times n,q)$ }

In this section, we consider bounds of $\chi_{d}(N\times n,q)$. Since the diameter of $M_{N\times n}(q)$ is $n$, $\chi_{d}(N\times n,q)= 1$ for any $d\ge n+1$. So we can assume $d\le n$.

\vskip.2cm

\begin{proposition}  For $N\geq n$, we have $\chi_{1}(N\times  n,q)= q^N$.
\end{proposition}
\p
When $d=1$,  it suffices to show that   $\chi_{d}(N\times n,q)\geq q^N$ since $\chi_{d}(N\times n,q)\leq \chi'_{d}(N,n,q)$.  Consider  a set $S=\{A=(a_{ij})\in \mathbb{F}^{N\times n}_{q}: a_{ij}=0,  \mbox{ for }1\leq i\leq N, 2\leq j\leq n \}$, in which  any two of them  have an edge, and so coloring them needs 
\[\left(\begin{array}{c}
N\\
0 \end{array} \right )+ \left(\begin{array}{c}
N\\
1 \end{array}\right )(q-1)+\ldots+\left(\begin{array}{c}
N\\
j \end{array}\right )(q-1)^j+\ldots+\left(\begin{array}{c}
N\\
N \end{array} \right )(q-1)^N=(1+q-1)^{N}.\]
\q


We give some examples about equidistant rank codes with maximal size in each case.

When $N=n=2$, $q=2$ and $d=2$. The code $C_1=\left( \left(\begin{array}{cc}
                                                                                                 1 & 0 \\
                                                                                                 0 & 1\end{array} \right ),
                                                                                                 \left(\begin{array}{cc}
                                                                                                 0 & 0 \\
                                                                                                 1 & 0\end{array} \right ),
                                                                                                 \left(\begin{array}{cc}
                                                                                                 0 & 1 \\
                                                                                                 0 & 0\end{array} \right ), \right.$
                                                                                                $  \left. \left(\begin{array}{cc}
                                                                                                 1 & 1 \\
                                                                                                 1 & 1\end{array} \right )
                                                                 \right )$ in matrix form.

Let
$\alpha$ be a primitive element of the Galois field $\mathbb{F}_{2^3}$ such that $\alpha^3=\alpha+1$. 
 We can construct an equidistant rank code 
 $C_2$ of length $n=2$ over $\mathbb{F}_{2^3}$ and rank distance $2$, where    
$C_2=\{(1,\alpha), (\alpha,1+\alpha), (1+\alpha+\alpha^2,1), (0,1+\alpha^2),
(1+\alpha,\alpha^2), (1+\alpha^2, 1+\alpha+\alpha^2),(\alpha+\alpha^2, \alpha+\alpha^2), (\alpha^2,0)\}$.  

We can also consider an example in \cite{sel}. Consider
an equidistance rank  code $C_3$ of length $n=3$ over $\mathbb{F}_{2^3}$ and rank distance $3$.   The code $C_3$ is  $\{(\alpha^2,0,0), (1,\alpha,\alpha^2), (0,1,\alpha), (1+\alpha^2, 1+\alpha,\alpha+\alpha^2), (\alpha+\alpha^2, \alpha+\alpha^2,1),
(1+\alpha,\alpha^2, 1+\alpha^2), (\alpha,1+\alpha+\alpha^2,1+\alpha), (1+\alpha+\alpha^2, 1+\alpha^2,1+\alpha+\alpha^2)\}$.

Therefore, we get the following proposition. 
\begin{proposition} We have 
\begin{equation} \label{eq:special} 
\chi_{n}(N\times n,2)= 2^N,  \mbox{  for pairs } (N,n)= (N,1),(2,2),(3,2),(3,3). \end{equation}
\end{proposition}

In addition, the size of an equidistant rank code over  $\mathbb{F}_{2^N}$with length $n$ and rank distance $n$    should be a power of 2 and so we propose  
an open problem as follow:
\begin{equation} \label{eq: conj}
\chi_{n}(N\times n,2)= 2^N \mbox{ for all }N\geq n. \end{equation}

If this statement holds, it also means that there is a largest equidistant rank code over $\mathbb{F}_{2^N}$ with length $n$, distance $n$, size $2^N$.

\vskip.2cm

\begin{proposition}
For $n\geq 3$, we have $\chi_{n-1}({n \choose 2}\times n ,q)\ge q^n-1$.
\end{proposition}
\p
By Proposition \ref{prop:5},
it is forward to have  $\chi_{n-1}({n \choose 2}\times n ,q)\ge q^n-1$.
\q

\vskip 0.2cm

 In the following we discuss upper bounds  of the exactly $d$-distance coloring of $\mathbb{F}^{N\times n}_q$. First we have 
 \begin{equation} \label{eq:nat}
 \chi_{d}(N\times  n,q)\leq \chi'_{d}(N\times  n,q)=q^{Nd}.
 \end{equation} 
 
 Furthermore, we use  the
 argument of linear forbidden rank distance codes.  An $[N\times n,k]_q$ linear rank code $\mathcal{M}$ is called an
$[N\times n,k,\overline{d}]_q$ linear forbidden rank distance code if $\mathcal{M}$ is also an $(N\times n, \overline{d})$ forbidden rank
distance code, i.e.,  $Rk(M)\neq d$ for any nonzero codeword $M$ of $\mathcal{M}$.
Given $N, n(n\leq N)$, and $d$, the maximum dimension of an
$[N\times n,k,\overline{d}]_q$ linear forbidden distance code is denoted by $k(N\times n,\overline{d})$. Similar to $d$-distance
coloring, we have

\begin{equation}
\chi_{d}(N\times  n,q)\leq q^{N(n-k)}.
\end{equation}
\noindent
\begin{equation}\label{eq:exactlc}
\chi_{d}(N\times  n,q)\leq q^{N(n-k(N\times n,\overline{d}))}.
\end{equation}

In the following lemma the lower bound on $k(N\times n,\overline{d})$ is given.

\vskip.2cm
\begin{lemma}  \label{lem:3.4}
We have
\begin{equation}
k(N\times n,\overline{d})\geq nN-\left\lceil \log_{q}[2+{n-1\choose d-1}(q^N-1)^{d-1}]\right\rceil.
\end{equation}
\end{lemma}
\p
First if
\begin{equation}
2+{n-1 \choose d-1}(q^N-1)^{d-1} \leq q^{mN},
\end{equation}
then there exists a  matrix $H$ over $\mathbb{F}_{q^N}$ of size $m\times n$ such that any column of $H$ is not equal to
the $\mathbb{F}_{q^N}$-linear combination of any other $d-1$ columns of $H$ by Lemma \ref{lem:3.3}.
Hence there is a $[N\times n, k,\overline{d}]$ linear forbidden rank distance code, where the dimension $k\geq n-m$. Take
\[Nm=\left\lceil \log_{q}[2+{n-1\choose d-1}(q^N-1)^{d-1}]\right\rceil. \]
Thus,
\begin{eqnarray*}
Nk(N, n,\overline{d}) &\geq & Ndim(\mathcal{V})\\
         & \geq & Nn-Nm \\
         & \geq & Nn-\left\lceil \log_{q}[2+{n-1\choose d-1}(q^N-1)^{d-1}]\right\rceil.
\end{eqnarray*}

This completes the proof.
\q

\vskip.2cm

\begin{theorem} We have
\begin{equation}
\chi_{d }(N\times n,q)\leq q^{\lceil \log_{q}[2+{n-1\choose d-1}(q^N-1)^{d-1}]\rceil}.
\end{equation}
\end{theorem}
\p
The theorem follows from Lemma \ref{lem:3.4} and Eq (\ref{eq:exactlc}).
\q

In the following table we list  some values of  two bounds from Eq(12) and Eq(\ref{eq:nat}). 

\begin{table}[h]
\setlength{\unitlength}{.68cm}
\begin{center}
\centering \caption{ \label{figfcsr}  \ \ \ \ Bounds Comparison }
\begin{tabular}{|c|c|c|c|c|c|} \hline
$N$ & $n$ & $d$ & $q$ & Bound $(12)$ & Bound (\ref{eq:nat})\\
\hline
$6$ & $4$ & $2$ & $2$ & $2^8$ & $2^{12}$\\
\hline
$6$ &  $4$ &  $3$ & $2$ & $2^{14}$ & $2^{18}$\\
\hline
$6$ &  $4$ &  $2$ & $3$ &  $3^7$  & $3^{12}$\\
\hline
$6$ &  $4$ &  $3$ & $3$ & $3^{13}$ & $3^{18}$\\
\hline
$5$ & $3$ & $2$ & $2$ & $2^{6}$ & $2^{10}$\\
\hline
$5$ & $3$ & $3$ & $3$ & $3^{10}$ & $3^{15}$\\
\hline
$10$ & $7$ & $4$ & $2$ & $2^{35}$ & $2^{40}$\\
\hline
$10$ & $7$ & $4$ & $3$ & $2^{33}$ & $2^{40}$\\
\hline
\end{tabular}
\end{center}
\end{table}



\vskip.2cm

\section*{Acknowledgements} This work is partially supported by National Natural Science Foundation of China
(Nos. 61373019, 11171097 and 61170289).


\vskip.4cm

\begin{thebibliography}{99}

\bibitem{Bhuyan} L.N. Bhuyan and D.P. Agrawal, Generalized hypercube and hyperbus structures for a computer network, IEEE Transactions on Computers, 33 (1984) 323-333.

\bibitem{del} P. Delsarte, Bilinear forms over a finite field, with applications to coding theory, 
Journal of Combinatorial theory A, 25 (1978) 226-241.

\bibitem{Tuvi} T. Etzion, N. Raviv, Equidistant codes in the Grassmannian, http://arxiv.org /abs/1308.6231.  

\bibitem{Fu} F.-W. Fu, S. Ling and C. Xing, New results on two hypercube coloring problem, Discrete Applied Mathematics,  161(2013)  2937-2945.

\bibitem{gab} E. M. Gabidulin, Theory of codes with maximum rank distance, Problems on Information Transmission,  21 (1) (1985) 1-12.

\bibitem{gab2008} E. M. Gabidulin, N. I. Pilipchuk, Error and erasure correcting algorithms for rank
codes, Designs, Codes and Cryptography, Springer Netherland, 49 (2008) 105-122.

\bibitem{gab2009} E. M. Gabidulin, M. Bossert, Algebraic codes in network coding, Probl. Inform. Trans.  45 (4) (2009) 3-17.

\bibitem{hanlu} Z. Han, M. Lu,  On two q-ary n-cube coloring problems,
http://arxiv.org/sub- mit/1386726.

\bibitem{hua} L.-K. Hua, A theorem on matrices over a field and its applications, Chinese Mathematical Society, 1 (2) (1951) 109-163.



















\bibitem{Jamison} R.E. Jamison, G.L. Matthews, Distance k colorings of Hamming graphs, Congr. Numer. 183 (2006) 193-202.

\bibitem{Kim} D.S. Kim, D.-Z. Du, P.M. Pardalos, A coloring problem on the n-cube, Discrete Appl. Math. 103 (2000) 307-311.

\bibitem{Klotz} W. Klotz, E. Sharifiyazdi, On the distance chromatic number of Hamming graphs, Adv. Appl. Discrete Math. 2 (2) (2008) 103-115.

\bibitem{Koetter} R. Koetter and F. R. Kschischang, Coding for errors and erasures
in random network coding, IEEE Transactions on Information Theory, 54 (2008)  3579-3591.

\bibitem{lidl}R. Lidl, H. Niederriter,  Finite Fields, Encyclopedia Math.
Appl. Vol. {\bf 20}, Addison-Wesley, Reading, 1983.

\bibitem{ling}S. Ling, L.J. Qu,  A note on linearized polynomials and the dimension of their
kernels, Finite Fields Appl., 18 (2012) 56-62.

\bibitem{Ngo} H.Q. Ngo, D.-Z. Du, R.L. Graham, New bounds on a hypercube coloring problem and linear codes, in: Proceedings of the International Conference on
Information Technology: Coding and Computing, ITCC'01, IEEE, 2001, 542-546.

\bibitem{NDG} H.Q. Ngo, D.-Z. Du and R.L. Graham, New bounds on a hypercube coloring problem, Inform. Process. Lett. 84 (2002) 265-269.

\bibitem{Ostergard} P.R.J. \"Osterg{\aa}rd, On a hypercube coloring problem, J. Combin. Theory Ser. A 108 (2004) 199-204.

\bibitem{Pavan} A. Pavan, P.-J. Wan, S.-R. Tong and D.H. Du, A new multihop lightwave network based on the generalized de-Bruijn graph, in: Proc. on IEEE
INFOCOM, 1996, 498-507.

\bibitem{sel} R. S. Selvaraj, J. Demamu, Equidistant rank metric codes: construction and properties, Communications in Information and Systems, 10 (3) (2010) 183-192.

\bibitem{Skupie} Z. Skupie\'n, BCH codes and distance multi- or fractional colorings in hypercubes asymptotically, Discrete Math. 307 (2007) 990-1000.

\bibitem{Sullivan} H. Sullivan and T. R. Bashkow, A large scale homogeneous full distributed parallel machine, I. Proceeding of 4th Annual Symposium on Computer Architecture, 1977, 105-117.

\bibitem{Wan} P.-J. Wan, Near-optimal conflict-free channel set assignments for an optical cluster-based hypercube network, J. Comb. Optim. 1 (1997) 179-186.


\end{thebibliography}
\end{document}